\documentclass[11pt]{article}
\usepackage[latin1]{inputenc}
\usepackage[colorlinks=true, bookmarksnumbered=true, bookmarksopen=true,
bookmarksopenlevel=3, pdfstartview=FitH, linkcolor=cyan, pdfmenubar=true,
pdftoolbar=true, bookmarks=true,citecolor=cyan, urlcolor=magenta,
filecolor=cyan,menucolor=black,plainpages=false,pdfpagelabels]{hyperref}
\usepackage[lmargin=2cm,tmargin=2cm,bmargin=2cm,rmargin=2cm]{geometry}
\usepackage{lineno,hyperref}
\usepackage{amsfonts}
\usepackage{amsmath,amssymb}
\usepackage{mathtools}
\usepackage{color}
\usepackage{amsmath,amssymb}
\usepackage{amsthm}
\usepackage{amsfonts}
\usepackage{tikz}
\usepackage{pstricks}
\usepackage{pstricks-add}
\usepackage{pst-all,pst-math,pst-plot}
\usepackage{textcomp}
\usepackage{graphicx,color}
\usepackage{graphicx}
\usepackage{fancyhdr,ifthen}

\usetikzlibrary{matrix}

\newtheorem{theorem}{Theorem}[section]
\newtheorem{lemma}[theorem]{Lemma}

\newtheorem{proposition}[theorem]{Proposition}

\newtheorem{remark}[theorem]{Remark}

\numberwithin{equation}{section}

\newcommand{\R}{\mathbb{R}}

\allowdisplaybreaks

\begin{document}
\title{Finite time blow-up for a multi-dimensional model of the Kiselev-Sarsam equation \\
\vspace{0.4cm}}
\author{
{Wanwan Zhang$^{1}$}{\thanks{{Corresponding author. }\newline
{E-mail address: zhangww@jxnu.edu.cn (Wanwan Zhang).}}} \\
\\
{\small $^{1}$ School of Mathematics and Statistics, Jiangxi Normal University,} \\
{\small {Nanchang, 330022, Jiangxi, P. R. China.}}
}
\date{}
\maketitle

\begin{abstract}
In this paper, we propose and investigate a multi-dimensional nonlocal active scalar equation of the form
\begin{eqnarray*}
\partial_t\rho+g\mathcal{R}_a\rho\cdot \nabla\rho= 0,~\rho(\cdot,0)=\rho_{0},
\end{eqnarray*}
where the transform $\mathcal{R}_a$ is defined by
\begin{eqnarray*}
\mathcal{R}_af(x)=\frac{\Gamma(\frac{n+1}{2})}{\pi^{\frac{n+1}{2}}}P.V.\int\limits_{\mathbb{R}^n}\Big(\frac{x-y}{|x-y|^{n+1}}-\frac{x-y}{(|x-y|^2+a^2)^{\frac{n+1}{2}}}\Big)f(y)dy.
\end{eqnarray*}
This model can be viewed as a natural generalization of the well-known Kiselev-Sasarm equation, which was introduced in \cite{[Kiselev-Sarsam]} as a one-dimensional model for the two-dimensional incompressible porous media equation.
We show the local well-posedness for this multi-dimensional model as well as the finite time gradient blow-up for a class of radial initial data.
\

\vspace{0.2cm}


\vspace{0.2cm}

\end{abstract}

\tableofcontents

\

\setlength{\headheight}{14.5pt} \setlength{\headwidth}{\textwidth} %
\pagestyle{fancy} \fancyhf{} \renewcommand{\headrulewidth}{0pt}
\chead{\ifthenelse{\isodd{\value{page}}}{Blow-up for multi-dimensional Kiselev-Sarsam model}{W. Zhang}} \rhead{\thepage}
\newpage

\section{Introduction and main results}
\hskip\parindent
The problem of finite time blow-up or global regularity for active scalar equations with nonlocal velocities
has attracted much attention during the last two decades.
We refer the readers to the paper \cite{[Kiselev]} for some classical examples of active scalar equations and related well-posedness results.

In this paper, we propose and  study the following multi-dimensional nonlocal active scalar equation
\begin{align}\label{M-KS}
\begin{cases}
\partial_t\rho+g\mathcal{R}_a\rho\cdot \nabla\rho= 0, &\, (x,t)\in \mathbb{R}^n\times\mathbb{R}_+,\\
\rho(x,0)=\rho_{0}(x),&\, x\in \mathbb{R}^n.
\end{cases}
\end{align}
Here $a,g>0$ are fixed constants, $n\geq2$ is the space dimension, and
the transform $\mathcal{R}_a=(\mathcal{R}^{(1)}_a,...,\mathcal{R}^{(n)}_a)$ is defined by
\begin{eqnarray}\label{definition-Ra}
\mathcal{R}_af(x)=P.V.\int\limits_{\mathbb{R}^n}\mathcal{K}_a(x-y)f(y)dy,
\end{eqnarray}
where
\begin{eqnarray}\label{Ra-kernel}
\mathcal{K}_a(x)
=\frac{\Gamma(\frac{n+1}{2})}{\pi^{\frac{n+1}{2}}}\Big(\frac{x}{|x|^{n+1}}-\frac{x}{(|x|^2+a^2)^{\frac{n+1}{2}}}\Big).
\end{eqnarray}
In the case when $n=1$, the equation \eqref{M-KS} is formally reduced to the famous Kiselev-Sarsam equation given by
\begin{align}\label{KS}
\begin{cases}
\partial_t\rho+g\mathcal{H}_a\rho\partial_x\rho= 0, &\, (x,t)\in \mathbb{R}\times\mathbb{R}_+,\\
\rho(x,0)=\rho_{0}(x),&\, x\in \mathbb{R},
\end{cases}
\end{align}
where the transform $\mathcal{H}_a$ is defined by
\begin{eqnarray}\label{definition-Ha}
\begin{split}
\mathcal{H}_af(x)
&=\frac{1}{\pi}P.V.\int\limits_{\mathbb{R}}\Big(\frac{1}{x-y}-\frac{x-y}{(x-y)^2+a^2}\Big)f(y)dy\\
&=\frac{1}{\pi}P.V.\int\limits_{\mathbb{R}}\frac{a^2f(y)}{(x-y)((x-y)^2+a^2)}dy.
\end{split}
\end{eqnarray}
We will construct the model \eqref{M-KS} in Section 2 by extending the one-dimensional transform $\mathcal{H}_a$ to the multi-dimensional transform $\mathcal{R}_a$.

The nonlocal active scalar transport equation \eqref{KS} was introduced by Kiselev and Sarsam in \cite{[Kiselev-Sarsam]} as a one-dimensional model analogy for the two-dimensional incompressible porous media (IPM) equation given by
\begin{align*}
\begin{cases}
\partial_t\rho+u\cdot\nabla\rho= 0, &\, (x,t)\in \mathbb{R}^2\times\mathbb{R}_+,\\
u=-\nabla P-(0,g\rho),~\nabla\cdot u=0,\\
\rho(x,0)=\rho_{0}(x),&\, x\in\mathbb{R}^2,
\end{cases}
\end{align*}
which models the transport of a scalar density $\rho(x,t)$ by an incompressible fluid velocity field $u(x,t)$ under the effects of Darcy's law and gravity. Here $P=P(x,t)$ is the scalar pressure, $g>0$ is the constant of gravitational acceleration, and $\rho_{0}=\rho_{0}(x)$ is the initial density.
The local well-posedness of the IPM equation has been established in various function spaces (see \cite{[Cordoba-Gancedo-Orive]} in H\"{o}lder spaces $C^{k,\beta}(\mathbb{R}^2)\cap L^2(\mathbb{R}^2)$ for $k\geq1$ and $0<\beta<1$, and \cite{[Constantin]} in Sobolev spaces $H^s(\mathbb{R}^2)$ for $s>2$).
In \cite{[Bianchini-Cordoba-Martninez]}, Bianchini, C\'{o}rdoba and Martn\'{i}nez-Zoroa showed the ill-posedness for IPM in the critical space $H^2(\mathbb{R}^2)$ and also proved the strong ill-posedness in $H^2(\mathbb{R}^2)$ for perturbations of a spectrally stable steady state.
While it seems that the question of finite-time singularity versus global regularity for IPM is challenging, there have been many works towards this problem. In \cite{[Kiselev-Yao]}, Kiselev and Yao constructed solutions to the IPM equation which exhibit infinite-in-time growth of derivatives provided global-in-time existence. In \cite{[Cordoba-Martninez-Zoroa]}, C\'{o}rdoba and Mart\'{i}nez-Zoroa showed the existence of a smooth compactly supported initial data and a compactly supported force in $L^\infty_tC^\infty_x$ that lead to finite-time blow-up in the forced IPM equation.
Furthermore, the Muskat equation, which models the evolution of the interface between two fluids of different densities evolving under the IPM equation, was studied by Zlato\v{s} in the setting of half-plane and was shown to exhibit finite-time singularity formation for smooth initial data \cite{[Zlatos-Local-regularity],[Zlatos-singularity-formation]}.
Very recently, in \cite{[Dembski]},  Dembski proved the finite-time singularity formation for Lipschitz continuous solutions of the IPM equation which vanish on the boundary of the domains arbitrarily close to the half-plane.
In \cite{[Kiselev-Sarsam]}, Kiselev and Sarsam gave detailed discussions on the derivation of the
one-dimensional model  \eqref{KS} from the IPM equation.
The authors also remarked how the 1-D IPM  equation \eqref{KS} parallels the well-known C\'{o}rdoba-C\'{o}rdoba-Fontelos (CCF) equation:
\begin{eqnarray}\label{CCF}
\partial_t\theta-\mathcal{H}\theta\partial_x\theta=0,~\theta(\cdot,0)=\theta_{0},
\end{eqnarray}
where $\mathcal{H}\theta$ is the Hilbert transform of $\theta$.

To motivate this paper, we first review some existing results on the CCF equation and its natural generalizations.
In \cite{[Cordoba-Cordoba-Fontelos05]}, C\'{o}rdoba-C\'{o}rdoba-Fontelos first showed the finite time singularity formation of solutions to \eqref{CCF} for a class of smooth even initial data.
Specifically, the blow-up proof in \cite{[Cordoba-Cordoba-Fontelos05]} is based on an ingenious inequality for the Hilbert transform $\mathcal{H}$: for any $-1<\delta<1$ and any even bounded smooth function $f$ defined on $\mathbb{R}$,
\begin{eqnarray}\label{CCF-inequality}
-\int\limits^\infty_0\frac{\mathcal{H} f(x) f'(x)}{x^{1+\delta}}dx
\geq
C_{\delta}\int\limits^\infty_0\frac{(f(x))^2}{x^{2+\delta}}dx,
\end{eqnarray}
where  $C_\delta>0$ is a constant depending only on $\delta$.  The proof of \eqref{CCF-inequality} in \cite{[Cordoba-Cordoba-Fontelos05]} is based on the Meillin transform and complex analysis.
In \cite{[Cordoba-Cordoba-Fontelos06]}, by establishing several new weighted nonlinear inequalities, C\'{o}rdoba et al. also proved the blow-up for positive, compactly supported initial data but that is not necessarily even. In \cite{[Castro-Cordoba]}, Castro and C\'{o}rdoba gave an simple short blow-up argument by considering even functions and by observing $\mathcal{H}(\mathcal{H}f f_{xx})(0)\geq0$.
Based on completely real variable arguments, Kiselev \cite{[Kiselev]} proved that more general inequality below
\begin{eqnarray}\label{nonlinear-inequality-H-Kiselev}
-\int\limits^1_0\frac{\mathcal{H} f(x) f'(x)(f(x))^{p-1}}{x^\delta}dx
\geq
C_{p,\delta}\int\limits^1_0\frac{(f(x))^{p+1}}{x^{1+\delta}}dx
\end{eqnarray}
holds true for any $p\geq1,\delta>0$ and any even bounded $C^1$ function $f$ with $f(0)=0$ and $f'\geq0$ on $(0,\infty)$.
In \cite{[Li-Rodrigo20]}, by using a pointwise inequality for the Hilbert transform, Li and Rodrigo gave several elementary proofs of the inequalities \eqref{CCF-inequality} and \eqref{nonlinear-inequality-H-Kiselev}. Silvestre and Vicol \cite{[Silvestre-Vicol]} provided four elegant blow-up proofs for the CCF equation \eqref{CCF}. In \cite{[Silvestre-Vicol]}, the authors also proved the finite time singularity of solutions to the fractionally transport velocity case ($\alpha$-CCF) given by
\begin{align}\label{CCF-generalized}
\partial_t\theta+(\Lambda^{-2+2\alpha}\partial_x\theta)\partial_x \theta= 0, ~\theta(\cdot,0)=\theta_{0},
\end{align}
for $0<\alpha<1$. When $\alpha=\frac12$, the equation \eqref{CCF-generalized} becomes the CCF equation \eqref{CCF}.

A multi-dimensional generalization of the CCF equation \eqref{CCF} given by
\begin{eqnarray}\label{M-CCF}
\partial_t\theta-\mathcal{R}\theta\cdot\nabla\theta=0,~\theta(\cdot,0)=\theta_{0},
\end{eqnarray}
was first considered in Balodis and C\'{o}rdaba \cite{[Balodis-Cordoba]}.
Here $\mathcal{R}\theta=(\mathcal{R}^{(1)}\theta,...,\mathcal{R}^{(n)}\theta)$ is the Riesz transform of $\theta$.
In \cite{[Balodis-Cordoba]}, the local well-posedness of solutions to \eqref{M-CCF} was established, and the authors also proved the finite-time blow-up of solutions to \eqref{M-CCF} by deriving and applying a multi-dimensional version of the integral inequality \eqref{CCF-inequality}: for any $c_n<\delta<1$ and a suitable smooth function $f$ with constant sign vanishing at the origin,
\begin{eqnarray*}
-\int\limits_{\R^n}\frac{(\mathcal{R} f(x)-\mathcal{R} f(0))\cdot\nabla f(x)}{|x|^{n+\delta}}dx
\geq
C_{\delta}\int\limits_{\R^n}\frac{(f(x))^2}{|x|^{n+1+\delta}}dx,
\end{eqnarray*}
where $0<c_n<1$ is a fixed constant.
When $n=2$, such result was also proved for a similar equation in \cite{[Dong-Li]} independently.
Later, the transport equation with  fractional velocity given by
\begin{align}\label{M-CCF-generalized}
\partial_t\theta+\Lambda^{-2+2\alpha}\nabla\theta\cdot \nabla\theta= 0, ~\theta(\cdot,0)=\theta_{0},
\end{align}
was also studied. Here the space dimension $n\geq2$ and $0<\alpha<1$.
The equation \eqref{M-CCF-generalized} with $\alpha=\frac12$ is reduced to the multi-dimensional CCF equation \eqref{M-CCF}. The local well-posedness of solutions to \eqref{M-CCF-generalized} in Sobolev spaces was established by Chae in \cite{[Chae]}. In \cite{[Dong]}, Dong  obtained the following weighted nonlinear inequality with full range $\alpha\in (0,1)$: for any $\delta\in(-2\alpha,2-2\alpha)$ and radial Schwartz function $f$,
\begin{eqnarray}\label{nonlinear-inequality}
\int\limits_{\R^n}\frac{\Lambda^{-2+2\alpha}\nabla f(x)\cdot\nabla f(x)}{|x|^{n+\delta}}dx
\geq
C_{n,\alpha,\delta}\int\limits_{\R^n}\frac{(f(0)-f(x))^2}{|x|^{n+2\alpha+\delta}}dx,
\end{eqnarray}
which was applied to prove the blow-up of smooth solutions to \eqref{M-CCF-generalized}  for any smooth, radially symmetric and nonnegative initial data with compact support and its positive maximum attained at the origin.
Under the radial and non-increasing assumption of $f$,  Li and Rodrigo \cite{[Li-Rodrigo20]} also proved the inequality \eqref{nonlinear-inequality} by deriving a pointwise inequality for the term $-\Lambda^{-2+2\alpha}\nabla f(x)\cdot\frac{x}{|x|}$ along with the use of the Hardy's inequality.
Motivated by \cite{[Silvestre-Vicol]}, Jiu and Zhang \cite{[Jiu-Zhang]} proved the finite time singularity of solutions to \eqref{M-CCF-generalized} for smooth initial data $\theta_0$ with $\displaystyle\sup_{x\in\R^n}\theta_0(x)>0$ via the De Giorgi iteration technique.

Finally, the finite-time blow-up problem of the fractionally dissipative equations of \eqref{CCF}, \eqref{CCF-generalized} and \eqref{M-CCF-generalized} in the supercritical scheme was also extensively investigated in the literature.
In summary, it was proved that certain solutions to the equation
\begin{align}\label{CCF-generalized-dissipation}
\begin{cases}
\partial_t\theta+(\Lambda^{-2+2\alpha}\partial_x\theta)\partial_x \theta+\Lambda^\gamma\theta= 0~ {\rm or}~\partial_t\theta+\Lambda^{-2+2\alpha}\nabla\theta\cdot\nabla \theta+\Lambda^\gamma\theta=0,\\
\rho(x,0)=\rho_{0}(x),
\end{cases}
\end{align}
develop finite time blow-up for when $\gamma\in(0,\alpha)$ for all $\alpha\in(0,1)$. We refer the readers to \cite{[Ferreira-Moitinho],[Kiselev],[Li-Rodrigo08],[Li-Rodrigo09],[Li-Rodrigo20],[Silvestre-Vicol],[Zhang]} and the
references therein for more details.
In the case when $\gamma\in[\alpha,2\alpha)$ for $\alpha\in(0,1)$, whether solutions of  \eqref{CCF-generalized-dissipation} with smooth initial data may blow up in finite time remains still open.

Now we proceed to review the results on the Kiselev-Sarsam equation \eqref{KS}. In \cite{[Kiselev-Sarsam]}, the authors proved the local well-posedness for the equation \eqref{KS} posed on the circle and adapted the arguments for the Hilbert transform in \cite{[Kiselev]} to show that for any $a,\delta>0$, it holds that
\begin{eqnarray}\label{IPM-inequality}
-\int\limits_0^{\frac\pi2}\frac{\mathcal{H}_af(x)f'(x)}{x^\delta}dx\geq C_{a,\delta}\int\limits_0^{\frac\pi2}\frac{(f(x))^2}{x^{1+\delta}}dx,
\end{eqnarray}
where $C_{a,\delta}$ is a universal constant depending only on $\delta$ and $a$, and $f$ is an even and nonnegative smooth function defined on $\mathbb{T}$  with $f(0)=0$ and $f'\geq0$ on $[0,\pi)$. As an application of the inequality \eqref{IPM-inequality}, Kiselev and Sarsam proved the finite time singularity of solutions to \eqref{KS} for a class of  smooth even initial data in the setting of the periodic circle. Recently, Liu and Zhang \cite{[Liu-Zhang]} established several weighted integral inequalities for the transform $\mathcal{H}_a$ in the setting of the real line. Based on these integral inequalities, the authors proved the finite time blow-up of solutions to \eqref{KS}.

It is then natural to consider the problem of finite time blow-up of the multi-dimensional Kiselev-Sarsam equation, which is currently absent form the literature, to our best of knowledge. The purpose of this paper is to introduce the multi-dimensional Kiselev-Sarsam equation \eqref{M-KS}, and prove the local well-posednesss for this model as well as the finite time blow-up for a class of radial initial data.



In Section 3, we establish our first result on the local well-posedness.
\begin{theorem}\label{local-well-posedness}
Let $n\geq2$ and $a,g>0$. For each $\rho_0\in H^s(\mathbb{R}^n)$ with $s>\frac n2+1$, there exists a $T=T(\|\rho_0\|_{H^s})>0$ such that \eqref{M-KS} admits a unique solution $\rho$ in $C([0,T); H^s(\mathbb{R}^n))\cap {\rm Lip}((0,T);H^{s-1}(\mathbb{R}^n))$.
\end{theorem}
Our main result shows that the family of smooth, compactly supported, radial and non-decreasing initial data with negative value at the origin to \eqref{M-KS} undergo finite time blow-up, whose proof is given in Section 4.
\begin{theorem}\label{singularity-formation-monotone}
Let $n\geq2$ and $a,g>0$.
Suppose $\rho_0\in C^\infty_c(\mathbb{R}^n)$ is radial and non-decreasing with $\rho_0(0)<0$.
Then the solution $\rho$ to \eqref{M-KS} with the initial data $\rho_0$ develops the gradient blow-up in finite time.
\end{theorem}

At the end of this section, some notations are introduced as follows.
For $p\in[1,\infty]$, we denote $L^p(\mathbb{R}^n)$ the standard Lebesgue space and its  norm by  $\|\cdot\|_{L^p(\mathbb{R}^n)}$.
For $s\geq0$, we use the notation $H^s(\mathbb{R}^n)$ to denote the nonhomogeneous Sobolev space of $s$ order,
whose endowed norm is denoted by $\|f\|_{H^s(\mathbb{R}^n)}=\|f\|_{L^2(\mathbb{R}^n)}+\|\Lambda^sf\|_{L^2(\mathbb{R}^n)}$, where the fractional Laplacian $\Lambda^s:=(-\Delta)^{\frac s2}$ is defined through the Fourier transform as
\begin{eqnarray*}
\widehat{(-\Delta)^{\frac s2}f}(\xi)=(2\pi|\xi|)^s\widehat{f}(\xi).
\end{eqnarray*}
$BMO(\mathbb{R}^n)$ denotes the space of functions  of bounded mean oscillation on $\mathbb{R}^n$ with the seminorm notation $\|\cdot\|_{BMO(\mathbb{R}^n)}$. For a sake of the convenience, the $L^{p}(\mathbb{R}^n)$-norm of a function $f$ is sometimes abbreviated as $\|f\|_{L^p}$, the $H^s(\mathbb{R}^n)$-norm as $\|f\|_{H^s}$ and the $BMO(\mathbb{R}^n)$-seminorm as $\|f\|_{BMO}$.
All norms of a function $f(x,t)$ depending on space and time variables will refer to the spatial norms.
Finally, the functions $\Gamma(\cdot)$ and $B(\cdot,\cdot)$ stand for the standard Gamma and Beta function, respectively.
Let $\mathbb{S}^{n-1}$ be the unit sphere in $\mathbb{R}^n$, i.e., $\mathbb{S}^{n-1}=\{x\in\mathbb{R}^n:|x|=1\}$
and $\omega_{n-1}$ be its surface area. We recall that
\begin{eqnarray*}
\omega_{n-1}=\frac{2\pi^{\frac{n}{2}}}{\Gamma(\frac{n}{2})}.
\end{eqnarray*}
Throughout this paper, we will use $C$ to denote a positive constant, whose value may change from line to line, and write $C_{n,a}$ or $C(n,a)$ to emphasize the dependence of a constant on $n$ and $a$.

The remaining part of this paper is organized as follows. In Section 2, we first introduce the multi-dimensional Kiselev-Sarsam equation \eqref{M-KS} by extending the one-dimensional transform $\mathcal{H}_a$ to the multi-dimensional transform $\mathcal{R}_a$, and then give some properties of $\mathcal{R}_a$.
Section 3 is devoted to the establishment of the local well-posedness of solutions to \eqref{M-KS}.
The proof of finite time blow-up of solutions is given in Section 4.
\section{Multi-dimensional extension of $\mathcal{H}_a$ and Properties of $\mathcal{R}_a$}
\hskip\parindent
The key ingredient of the construction of the multi-dimensional Kiselev-Sarsam equation \eqref{M-KS} is to extend the one-dimensional transform $\mathcal{H}_a$ to the multi-dimensional transform $\mathcal{R}_a$. For every $a>0$, by \eqref{definition-Ha}, we know that the transform $\mathcal{H}_a$ is a convolution operator, and note that the corresponding kernel $K_a(x):=\frac{1}{\pi x}\frac{a^2}{x^2+a^2}$ can be represented as the Hilbert kernel $\frac{1}{\pi x}$ minus the one-dimensional conjugate Poisson kernel
\begin{align*}
Q_a(x):=\frac{1}{\pi}\frac{x}{x^2+a^2}.
\end{align*}
As a natural generalization of $\mathcal{H}_a$, the kernel $\mathcal{K}_a(x)$ of the convolution operator $\mathcal{R}_a$ should be defined as \eqref{Ra-kernel}, which is exactly the difference between the Riesz kernel $\frac{\Gamma(\frac{n+1}{2})}{\pi^{\frac{n+1}{2}}}\frac{x}{|x|^{n+1}}$  and the multi-dimensional conjugate Poisson kernel
\begin{align*}
\mathcal{Q}_a(x):=\frac{\Gamma(\frac{n+1}{2})}{\pi^{\frac{n+1}{2}}}\frac{x}{(|x|^2+a^2)^{\frac{n+1}{2}}}.
\end{align*}

On the other hand, being a convolution operator, the transform $\mathcal{H}_a$ is also a Fourier multiplier operator on the real line with the symbol
\begin{align}\label{1d-symbol}
\widehat{K_a}(\eta)=-i~{\rm sgn}(\eta)(1-e^{-2\pi a|\eta|})=-i~{\rm sgn}(\eta)+i~{\rm sgn}(\eta)e^{-2\pi a|\eta|}.
\end{align}
We note that the Fourier transform of the multi-dimensional conjugate Poisson kernel $\mathcal{Q}_a(x)$ is
\begin{align*}
\widehat{\mathcal{Q}_a}(\xi)=-\frac{i\xi}{|\xi|}e^{-2\pi a|\xi|}
\end{align*}
(see, e.g., Exercise 5.1.8 in \cite{[Grafakos]}). It follows that the transform $\mathcal{R}_a$ given by \eqref{definition-Ra} is also a Fourier multiplier operator on the whole space $\mathbb{R}^n$ with the symbol
\begin{align*}
\widehat{\mathcal{K}_a}(\xi)=-\frac{i\xi}{|\xi|}-\Big(-\frac{i\xi}{|\xi|}e^{-2\pi a|\xi|}\Big)
=-\frac{i\xi}{|\xi|}(1-e^{-2\pi a|\xi|}),
\end{align*}
which is exactly the multi-dimensional version of the one-dimensional symbol \eqref{1d-symbol}.

Altogether, the transform $\mathcal{R}_a$ defined by \eqref{definition-Ra} may indeed be a reasonable extension of the one-dimensional transform $\mathcal{H}_a$. Thus, the multi-dimensional model \eqref{M-KS} can be viewed as a natural generalization of the Kiselev-Sarsam equation \eqref{KS}.

In addition, one can think of $\mathcal{R}_a$ as an operator that interpolates between the trivial zero operator and the prototypical singular integral operator: the Riesz transform $\mathcal{R}$ with kernel $\frac{\Gamma(\frac{n+1}{2})}{\pi^{\frac{n+1}{2}}}\frac{x}{|x|^{n+1}}$.
This can be seen in two ways. First, the kernel $\mathcal{K}_a$ converges pointwise to the Riesz kernel as $a\rightarrow\infty$ while it instead converges pointwise to zero when taking $a\rightarrow0$.

Second, the symbol $\widehat{\mathcal{K}_a}(\xi)$ of the Fourier multiplier operator $\mathcal{R}_a$ converges pointwise to the symbol of the Riesz transform as $a\rightarrow\infty$, while it instead converges pointwise to zero as $a\rightarrow0$. Therefore, by the dominated convergence theorem, it holds that
\begin{align*}
\|(\mathcal{R}_a-\mathcal{R})f\|_{L^2}
=\Big\|\frac{i\xi}{|\xi|}e^{-2\pi a|\xi|}\widehat{f}\Big\|_{L^2}
=\Big\|e^{-2\pi a|\xi|}\widehat{f}\Big\|_{L^2}\rightarrow0~{\rm as} ~a\rightarrow\infty
\end{align*}
and
\begin{align*}
\|\mathcal{R}_af\|_{L^2}
=\Big\|-\frac{i\xi}{|\xi|}\Big(1-e^{-2\pi a|\xi|}\Big)\widehat{f}\Big\|_{L^2}
=\Big\|\Big(1-e^{-2\pi a|\xi|}\Big)\widehat{f}\Big\|_{L^2}\rightarrow0~{\rm as} ~a\rightarrow0.
\end{align*}
In words, we have that $\mathcal{R}_a$ converges to $\mathcal{R}$ as $a\rightarrow\infty$ while instead converging to zero as $a\rightarrow0$, both with respect to the $L^2$ strong operator topology.

Finally, we show some bounded properties of the transform $\mathcal{R}_a$, which will be needed later. Since
\begin{align}\label{L2-boundness-Ra}
\|\mathcal{R}_af\|_{L^2}
=\Big\|\Big(1-e^{-2\pi a|\xi|}\Big)\widehat{f}\Big\|_{L^2}\leq\|f\|_{L^2},
\end{align}
and then
\begin{align}\label{homogeneous-Hs-boundness-Ra}
\|\Lambda^s\mathcal{R}_af\|_{L^2}
=\|\mathcal{R}_a\Lambda^sf\|_{L^2}\leq\|\Lambda^sf\|_{L^2}.
\end{align}
These mean that, for any $s\geq0$,
\begin{align}\label{Hs-boundness-Ra}
\|\mathcal{R}_af\|_{H^s}\leq\|f\|_{H^s},
\end{align}
which together with the continuous embedding $H^\lambda(\mathbb{R}^n)\hookrightarrow L^\infty(\mathbb{R}^n)$ for $\lambda>\frac n2$ implies that
\begin{eqnarray}\label{control}
\begin{split}
\|\nabla f\|_{L^\infty}+\|{\rm div} \mathcal{R}_af\|_{L^\infty}+\|\nabla \mathcal{R}_af\|_{L^\infty}
&\leq C_{n,s}(\|\nabla f\|_{H^{s-1}}+\|\partial_k \mathcal{R}^{(k)}_af\|_{H^{s-1}}+\|\partial_j \mathcal{R}^{(k)}_af\|_{H^{s-1}})\\
&\leq C_{n,s}(\|\nabla f\|_{H^{s-1}}+\|\partial_k f\|_{H^{s-1}}+\|\partial_j f\|_{H^{s-1}})\\
&\leq C_{n,s}\|f\|_{H^s},
\end{split}
\end{eqnarray}
for any $f\in H^s(\mathbb{R}^n)$ with $s>\frac{n}{2}+1$.

Also, the transform $\mathcal{R}_a$ satisfies the assumptions of Calder\'{o}n-Zygmund theory, with being a bounded linear operator on $L^p$ for any $p\in(1,\infty)$. In addition, the transform $\mathcal{R}_a$ maps $L^\infty(\mathbb{R}^n)$ to $BMO(\mathbb{R}^n)$, that is
\begin{align}\label{BMO-L-infinity}
\|\mathcal{R}_af\|_{L^\infty}\leq C_{n,a}\|f\|_{BMO}.
\end{align}
We refer the readers to \cite{[Abels]} and \cite{[Stein1970]} for the details of the proof of these properties for the singular integral operator $\mathcal{R}_a$.
\section{Local well-posedness}
\hskip\parindent In this section, we will present the lemmas required to prove Theorem \ref{local-well-posedness}. We first prove the uniqueness of solutions to \eqref{M-KS}.
\begin{lemma}\label{uniqueness}
Fix $n\geq2$ and $a,g>0$. Suppose $\rho_1,\rho_2$ are solutions in $C([0,T); H^s(\mathbb{R}^n))$ to \eqref{M-KS} with respect to the same initial data $\rho_0\in H^s(\mathbb{R}^n)$  for some $s>\frac{n}{2}+1$.  Then $\rho_1=\rho_2$ on $[0,T)$.
\end{lemma}
\textbf{Proof.} All computations and estimates below hold on the time interval $[0,T)$.
Denoting $\widetilde{\rho}:=\rho_1-\rho_2$, by \eqref{M-KS}, we have that
\begin{eqnarray*}
\frac{1}{g}\partial_t\widetilde{\rho}=-\mathcal{R}_a\rho_1\cdot\nabla\rho_1+\mathcal{R}_a\rho_2\cdot\nabla\rho_2=-\mathcal{R}_a\rho_1\cdot\nabla\widetilde{\rho}-\mathcal{R}_a\widetilde{\rho}\cdot\nabla\rho_2.
\end{eqnarray*}
It follows that
\begin{align*}
\frac{1}{2g}\frac{d}{dt}\|\widetilde{\rho}\|^2_{L^2}=-\int\limits_{\mathbb{R}^n}\widetilde{\rho}\mathcal{R}_a\rho_1\cdot\nabla\widetilde{\rho} dx-\int\limits_{\mathbb{R}^n}\widetilde{\rho} \mathcal{R}_a\widetilde{\rho}\cdot\nabla\rho_2dx.
\end{align*}
We observe that
\begin{align*}
-\int\limits_{\mathbb{R}^n}\widetilde{\rho}\mathcal{R}_a\rho_1\cdot\nabla\widetilde{\rho} dx
=\frac12\int\limits_{\mathbb{R}^n}\widetilde{\rho}^2{\rm div} \mathcal{R}_a\rho_1dx\leq\|{\rm div} \mathcal{R}_a\rho_1\|_{L^\infty}\|\widetilde{\rho}\|^2_{L^2}
\leq C_{n,s}\|\rho_1\|_{H^s}\|\widetilde{\rho}\|^2_{L^2},
\end{align*}
where the final inequality holds for any $s>\frac{n}{2}+1$ by \eqref{control}, and by \eqref{L2-boundness-Ra} and \eqref{control}
\begin{align*}
\Big|\int\limits_{\mathbb{R}^n}\widetilde{\rho}\mathcal{R}_a\widetilde{\rho}\cdot\nabla\rho_2dx\Big|
\leq\|\widetilde{\rho}\|_{L^2}\|\mathcal{R}_a\widetilde{\rho}\|_{L^2}\|\nabla\rho_2\|_{L^\infty}
\leq C_{n,s}\|\rho_2\|_{H^s}\|\widetilde{\rho}\|^2_{L^2}.
\end{align*}
Altogether, we obtain that
\begin{eqnarray*}
\frac{d}{dt}\|\widetilde{\rho}\|^2_{L^2}
\leq gC_{n,s}(\|\rho_1\|_{H^s}+\|\rho_2\|_{H^s})\|\widetilde{\rho}\|^2_{L^2}.
\end{eqnarray*}
Gr\"{o}nwall's inequality along with $\widetilde{\rho}(x,0)=0$ finishes the proof of the uniqueness of solutions.\hfill\hfill$\square$\vskip12pt
We next establish a-priori estimates on the $L^2$ norm of a solution.
\begin{lemma}\label{L-2}
Fix $n\geq2$ and $a,g>0$.  Suppose $\rho$ is a solution to \eqref{M-KS} in $C([0,T); H^s(\mathbb{R}^n))$ for some $s>\frac{n}{2}+1$. It then holds that, for any $t\in[0,T)$,
\begin{align*}
\frac{1}{g}\frac{d}{dt}\|\rho\|^2_{L^2}\leq C_{n,s}\|\rho\|_{H^s}\|\rho\|^2_{L^2}.
\end{align*}
\end{lemma}
\textbf{Proof.} We have by \eqref{M-KS} and \eqref{control} that
\begin{eqnarray}\label{L-2-estimate}
\begin{split}
\frac{1}{g}\frac{d}{dt}\|\rho\|^2_{L^2}
&=-2\int\limits_{\mathbb{R}^n}\rho \mathcal{R}_a\rho\cdot\nabla\rho dx
=\int\limits_{\mathbb{R}^n}\rho^2{\rm div} \mathcal{R}_a\rho dx\\
&\leq\|{\rm div} \mathcal{R}_a\rho\|_{L^\infty}\|\rho\|^2_{L^2}\leq C_{n,s}\|\rho\|_{H^s}\|\rho\|^2_{L^2}
\end{split}
\end{eqnarray}
on the time interval $[0,T)$.\hfill\hfill$\square$\vskip12pt

We proceed to bound the $\dot{H}^s$ seminorm of a solution. To do so, we make use of the following Kato-Ponce commutator estimates, whose proof can be found in \cite{[Ju]}.
\begin{lemma}\label{commutator-estimates}
Let $s>0$ and let $p,p_1,p_4\in(1,\infty)$, $p_2,p_3\in(1,\infty]$ such that $\frac{1}{p}=\frac{1}{p_1}+\frac{1}{p_2}=\frac{1}{p_3}+\frac{1}{p_4}$.
For $f,g\in \mathcal{S}(\mathbb{R}^n)$, there exists a constant $C>0$ depending only on $n,s,p,p_1$ and $p_3$ such that
\begin{eqnarray*}
\|\Lambda^{s}(fg)-f\Lambda^{s}g\|_{L^p}\leq C \Big(\|\Lambda^{s}f\|_{L^{p_1}}\|g\|_{L^{p_2}}+\|\nabla f\|_{L^{p_3}}\|\|\Lambda^{s-1}g\|_{L^{p_4}}\Big).
\end{eqnarray*}
\end{lemma}
Then we have
\begin{lemma}\label{derivative-estimates}
Fix $n\geq2$ and $a,g>0$. Suppose $\rho$ is a solution to \eqref{M-KS} in $C([0,T); H^s(\mathbb{R}^n))$ for some $s>\frac{n}{2}+1$. It then holds that
\begin{align*}
\frac{1}{2g}\frac{d}{dt}\|\rho\|^2_{\dot{H}^s}\leq C_{n,a,s}\|\rho\|_{H^s}\|\rho\|^2_{\dot{H}^s}
\end{align*}
on $[0,T)$.
\end{lemma}
\textbf{Proof.} All computations and estimates below hold on the time interval $[0,T)$.
We observe that
\begin{align*}
\frac{1}{2g}\frac{d}{dt}\|\Lambda^s\rho\|^2_{L^2}
&=-\int\limits_{\mathbb{R}^n}\Lambda^s\rho\Lambda^s(\mathcal{R}_a\rho\cdot\nabla\rho)dx\\
&=-\int\limits_{\mathbb{R}^n}\Lambda^s\rho\Big(\Lambda^s(\mathcal{R}_a\rho\cdot\nabla\rho)-\mathcal{R}_a\rho\cdot\Lambda^s\nabla\rho\Big) dx
+\frac 12\int\limits_{\mathbb{R}^n}(\Lambda^s\rho)^2{\rm div} \mathcal{R}_a\rho dx\\
&\leq \|\Lambda^s\rho\|_{L^2}\|\Lambda^s(\mathcal{R}_a\rho\cdot\nabla\rho)-\mathcal{R}_a\rho\cdot\Lambda^s\nabla\rho\|_{L^2}
+\|{\rm div} \mathcal{R}_a\rho\|_{L^\infty}\|\Lambda^s\rho\|^2_{L^2}.
\end{align*}
By Lemma \ref{commutator-estimates} and \eqref{homogeneous-Hs-boundness-Ra}, we can estimate the commutator as
\begin{align*}
\|\Lambda^s(\mathcal{R}_a\rho\cdot\nabla\rho)-\mathcal{R}_a\rho\cdot\Lambda^s\nabla\rho\|_{L^2}
&\leq C_{n,s}\Big(\|\Lambda^s\mathcal{R}_a\rho\|_{L^2}\|\nabla\rho\|_{L^\infty}
+\|\Lambda^{s-1}\nabla\rho\|_{L^2}\|\nabla \mathcal{R}_a\rho\|_{L^\infty}\Big)\\
&\leq C_{n,s}\Big(\|\nabla\rho\|_{L^\infty}
+\|\nabla \mathcal{R}_a\rho\|_{L^\infty}\Big)\|\Lambda^s\rho\|_{L^2}.
\end{align*}
Therefore, we obtain that
\begin{eqnarray}\label{High-order-derivative-estimate}
\frac{1}{2g}\frac{d}{dt}\|\rho\|^2_{\dot{H}^s}
\leq C_{n,s}\Big(\|\nabla\rho\|_{L^\infty}
+\|\nabla \mathcal{R}_a\rho\|_{L^\infty}+\|{\rm div} \mathcal{R}_a\rho\|_{L^\infty}\Big)\|\rho\|^2_{\dot{H}^s}.
\end{eqnarray}
which along with \eqref{control} yields that
\begin{eqnarray*}
\frac{1}{2g}\frac{d}{dt}\|\rho\|^2_{\dot{H}^s}
\leq C_{n,s}\|\rho\|_{H^s}\|\rho\|^2_{\dot{H}^s},
\end{eqnarray*}
which completes the proof Lemma \ref{derivative-estimates}. \hfill\hfill$\square$\vskip12pt
Lastly, it remains to bound the ${\rm Lip}((0,T);H^{s-1}(\mathbb{R}^n))$  norm of a solution.
\begin{lemma}\label{Lipschitz}
Fix $n\geq2$ and $a,g>0$.
Suppose $\rho$ is a solution to \eqref{M-KS} in $C([0,T); H^s(\mathbb{R}^n))$ for some $s>\frac{n}{2}+1$.
Then,  we have
\begin{align*}
\|\rho\|_{{\rm Lip}((0,T);H^{s-1}(\mathbb{R}^n))}\leq gC_{n,s}\|\rho\|^2_{L^\infty((0,T); H^s)}.
\end{align*}
\end{lemma}
\textbf{Proof.} By \eqref{M-KS}, the fact that $H^{s-1}(\mathbb{R}^n)$ is an algebra, and \eqref{Hs-boundness-Ra}, we have
\begin{align*}
\|\partial_t\rho\|_{H^{s-1}}
=g\|\mathcal{R}_a\rho\cdot \nabla\rho\|_{H^{s-1}}
\leq gC_{n,s}\|\mathcal{R}_a\rho\|_{H^{s-1}}\|\nabla\rho\|_{H^{s-1}}
\leq gC_{n,s}\|\rho\|^2_{H^s}.
\end{align*}
Therefore, for all $0<t_1<t_2<T$,
\begin{align*}
\|\rho(t_2)-\rho(t_1)\|_{H^{s-1}}\leq\int\limits^{t_2}_{t_1}\|\partial_t\rho(t)\|_{H^{s-1}}dt\leq gC_{n,s}(t_2-t_1)\|\rho\|^2_{L^\infty((0,T); H^s)},
\end{align*}
which concludes the proof of Lemma \ref{Lipschitz}. \hfill\hfill$\square$\vskip12pt

\textbf{Proof of Theorem \ref{local-well-posedness}.}
Collecting Lemmas \ref{L-2} and \ref{derivative-estimates}, one can get
\begin{eqnarray*}
\frac{d}{dt}\|\rho\|_{H^s}
\leq gC_{n,s}\|\rho\|^2_{H^s},
\end{eqnarray*}
which implies that
\begin{align*}
\|\rho(t)\|_{H^s}\leq\frac{\|\rho_0\|_{H^s}}{1-gC_{n,s}\|\rho_0\|_{H^s}t},~t\in\Big[0,\frac{1}{gC_{n,s}\|\rho_0\|_{H^s}}\Big).
\end{align*}
This provides us with a fundamental a priori estimate for \eqref{M-KS} in the $H^s$ norm,
\begin{align*}
\|\rho\|_{L^\infty((0,T); H^s)}\leq2\|\rho_0\|_{H^s}, ~{\rm where~}T:=\frac{1}{2gC_{n,s}\|\rho_0\|_{H^s}},
\end{align*}
which along with Lemma \ref{Lipschitz} yields that
\begin{align*}
\|\rho\|_{{\rm Lip}((0,T);H^{s-1}(\mathbb{R}^n))}\leq 4gC_{n,s}\|\rho_0\|^2_{H^s}.
\end{align*}
We thus obtain local-in-time a priori estimates in $L^\infty((0,T); H^s(\mathbb{R}^n))\cap {\rm Lip}((0,T);H^{s-1}(\mathbb{R}^n))$.
Then we can establish the local existence of a solution to \eqref{M-KS} by the standard argument of approximation by mollification.
Specifically, one needs to work with the regularized system
\begin{align*}
\begin{cases}
\partial_t\rho^\epsilon+gJ_\epsilon(\mathcal{R}_aJ_\epsilon\rho^\epsilon\cdot\nabla J_\epsilon\rho^\epsilon)=0, &\, (x,t)\in \mathbb{R}^n\times\mathbb{R}_+,\\
\rho^\epsilon(x,0)=\rho_{0}(x),&\, x\in \mathbb{R}^n,
\end{cases}
\end{align*}
where $J_\epsilon$ is the standard mollifier. For a sake of conciseness, we leave the interested reader to check the details, which can be consulted in \cite{[Chae],[Maida-Bertozzi]}.
Finally, Lemma \ref{uniqueness} ensures the uniqueness of the solution, and hence Theorem \ref{local-well-posedness} holds true for any fixed choice of $n\geq2$ and $a,g>0$.
\section{Finite time blow-up of solutions}
\hskip\parindent In this section, we prove the finite time blow-up of smooth solutions to \eqref{M-KS} for a class of radial smooth initial data.
\subsection{B-K-M type criterion and some properties of the solution}
Now that we have established local well-posedness, we assert the following Beale-Kato-Majda type criterion for \eqref{M-KS}. Before that, we recall a limiting Sobolev inequality needed later, which was proved in \cite{[Kozono-Taniuchi]}.
\begin{lemma}\label{limiting-Sobolev-inequality}
Let $1<p<\infty$ and let $s>\frac{n}{p}$. There is a constant $C=C(n,p,s)$ such that the estimate
\begin{align*}
\|f\|_{L^\infty}\leq C(1+\|f\|_{BMO})(1+\log^+\|f\|_{W^{s,p}})
\end{align*}
holds for all $f\in W^{s,p}(\mathbb{R}^n)$.
\end{lemma}
Then we have
\begin{proposition}\label{Beale-Kato-Majda-type}
Fix $n\geq2$ and $a,g>0$. Suppose $\rho$ is a solution to \eqref{M-KS} in $C([0,T_\ast);H^s(\mathbb{R}^n))$ corresponding to an initial data $\rho_0\in H^s(\mathbb{R}^n)$ for some $s>\frac{n}{2}+1$. If $0<T_\ast<\infty$ is the first blow-up time such that $\rho$ cannot be continued in $C([0,T_\ast); H^s(\mathbb{R}^n))$, then we must have that
\begin{eqnarray*}
\displaystyle\limsup_{t\rightarrow T_\ast}\|\rho(t)\|_{H^s}=\infty~{\rm if~and~only~if}~\lim_{t\rightarrow T_\ast}\int\limits_0^{t}\|\nabla\rho(\tau)\|_{L^\infty}d\tau=\infty.
\end{eqnarray*}
\end{proposition}
\textbf{Proof.} By \eqref{High-order-derivative-estimate} and the first inequality in \eqref{L-2-estimate}, we have
\begin{eqnarray}\label{Hs-estimate}
\frac{1}{2}\frac{d}{dt}\|\rho\|^2_{H^s}\leq g C_{n,s}(\|\nabla\rho\|_{L^\infty}+\|\nabla \mathcal{R}_a\rho\|_{L^\infty}+\|{\rm div} \mathcal{R}_a\rho\|_{L^\infty})\|\rho\|^2_{H^s}.
\end{eqnarray}
By Lemma \ref{limiting-Sobolev-inequality} and the boundedness properties \eqref{Hs-boundness-Ra} and \eqref{BMO-L-infinity}, we have that, for $s>\frac{n}{2}+1$,
\begin{eqnarray*}
\begin{split}
\|\partial_k \mathcal{R}^{(j)}_a\rho\|_{L^\infty}
&\leq C_{n,s}(1+\|\mathcal{R}^{(j)}_a\partial_k\rho\|_{BMO})(1+\log^+\|\mathcal{R}^{(j)}_a\partial_k\rho\|_{H^{s-1}})\\
&\leq C_{n,a,s}(1+\|\nabla\rho\|_{L^\infty})\log(e+\|\rho\|_{H^{s}})
\end{split}
\end{eqnarray*}
and
\begin{eqnarray*}
\begin{split}
\|\partial_j \mathcal{R}^{(j)}_a\rho\|_{L^\infty}
&\leq C_{n,s}(1+\|\mathcal{R}^{(j)}_a\partial_j\rho\|_{BMO})(1+\log^+\|\mathcal{R}^{(j)}_a\partial_j\rho\|_{H^{s-1}})\\
&\leq C_{n,a,s}(1+\|\nabla\rho\|_{L^\infty})\log(e+\|\rho\|_{H^{s}}).
\end{split}
\end{eqnarray*}
Substituting these logarithmic-type estimates into \eqref{Hs-estimate}, we obtain that
\begin{eqnarray*}
\frac{1}{2}\frac{d}{dt}\|\rho\|^2_{H^s}\leq gC_{n,a,s}(1+\|\nabla\rho\|_{L^\infty})\log(e+\|\rho\|_{H^{s}})\|\rho\|^2_{H^s},
\end{eqnarray*}
which implies that
\begin{align*}
\frac{d}{dt}\log(e+\|\rho(t)\|_{H^s})\leq gC_{n,a,s}(1+\|\nabla\rho\|_{L^\infty})\log(e+\|\rho\|_{H^{s}}).
\end{align*}
Therefore, it follows from Gr\"{o}nwall's inequality that
\begin{align*}
\|\rho(t)\|_{H^s}\leq  (e+\|\rho_0\|_{H^s})^{\exp\{Cg\int_0^t(1+\|\nabla\rho(\tau)\|_{L^\infty})d\tau\}}.
\end{align*}
Conversely, for $s>\frac{n}{2}+1$, by \eqref{control}, we have that
\begin{align*}
\int\limits_0^t\|\nabla\rho(\tau)\|_{L^\infty}d\tau\leq Ct\sup_{0\leq\tau\leq t}\|\rho(\tau)\|_{H^s}.
\end{align*}
The proof of Proposition \ref{Beale-Kato-Majda-type} is now finished. \hfill\hfill$\square$\vskip12pt
We proceed to show the property of solutions that the radial symmetry and non-decreasing monotonicity of the initial data can be preserved by the solution to \eqref{M-KS}.
\begin{lemma}\label{radial-property-preserved}
Fix $n\geq2$ and $a,g>0$. If $\rho$ is a smooth solution to \eqref{M-KS} with a radial and non-decreasing initial data $\rho_0$, then  $\rho(x,t)$ is also radial and non-decreasing in its life span.
\end{lemma}
\textbf{Proof}. Let $\mathbf{O}\in\mathbb{R}^{n\times n}$ be any orthogonal matrix. By the uniqueness of solutions to \eqref{M-KS} and  the radial property of $\rho_0$, it suffices to show that the function $\rho_{\mathbf{O}}(x,t):=\rho(\mathbf{O}x,t)$
is also a solution to \eqref{M-KS} with the initial data $\rho_0(\mathbf{O}x)$. Indeed, standard computations give that
\begin{eqnarray*}
(\partial_t\rho_{\mathbf{O}})(x,t)=(\partial_t\rho)(\mathbf{O}x,t),~(\nabla_x\rho_{\mathbf{O}})(x,t)={\mathbf{O}}^T(\nabla\rho)(\mathbf{O}x,t).
\end{eqnarray*}
By \eqref{definition-Ra} and \eqref{Ra-kernel}, we can derive that
\begin{align*}
\mathcal{R}_a\rho_{\mathbf{O}}(x,t)
&=\frac{\Gamma(\frac{n+1}{2})}{\pi^{\frac{n+1}{2}}}P.V.\int\limits_{\mathbb{R}^{n}}\Big(\frac{x-y}{|x-y|^{n+1}}-\frac{x-y}{(|x-y|^2+a^2)^{\frac{n+1}{2}}}\Big)\rho(\mathbf{O}y, t)dy\nonumber\\
&=\frac{\Gamma(\frac{n+1}{2})}{\pi^{\frac{n+1}{2}}}P.V.\int\limits_{\mathbb{R}^{n}}\Big(\frac{x-{\mathbf{O}}^{-1}z}{|x-{\mathbf{O}}^{-1}z|^{n+1}}-\frac{x-{\mathbf{O}}^{-1}z}{(|x-{\mathbf{O}}^{-1}z|^2+a^2)^{\frac{n+1}{2}}}\Big)\rho(z, t)dz \nonumber \\
&=\frac{\Gamma(\frac{n+1}{2})}{\pi^{\frac{n+1}{2}}} {\mathbf{O}}^{-1}P.V.\int\limits_{\mathbb{R}^{n}}\Big(\frac{\mathbf{O} x-z}{|\mathbf{O} x-z|^{n+1}}-\frac{\mathbf{O}x-z}{(|\mathbf{O}x-z|^2+a^2)^{\frac{n+1}{2}}}\Big)\rho(z,t)dz \nonumber \\
&={\mathbf{O}}^{-1} \mathcal{R}_a\rho(\mathbf{O} x,t).
\end{align*}
Thus, we obtain
\begin{align*}
(\partial_t\rho_{\mathbf{O}}+g\mathcal{R}_a\rho_{\mathbf{O}}\cdot\nabla\rho_{\mathbf{O}})(x,t)
&=(\partial_t\rho)(\mathbf{O}x,t)+g{\mathbf{O}}^{-1} \mathcal{R}_a\rho(\mathbf{O} x,t)\cdot {\mathbf{O}}^T(\nabla\rho)(\mathbf{O}x,t)\\
&=(\partial_t\rho)(\mathbf{O}x,t)+g\mathcal{R}_a\rho(\mathbf{O} x,t)\cdot \mathbf{O}{\mathbf{O}}^T(\nabla\rho)(\mathbf{O}x,t)\\
&=(\partial_t\rho+g\mathcal{R}_a\rho\cdot\nabla\rho)(\mathbf{O}x,t)=0.
\end{align*}
For any radially symmetric solution $\rho(x,t)=\rho(|x|,t)$ to \eqref{M-KS} with a radial and non-decreasing initial data $\rho_0(x)=\rho_0(|x|)$, by using polar coordinates, the equation \eqref{M-KS} is then reduced to
\begin{align}\label{radial-equation}
\partial_t\rho(r,t)+g\widetilde{\mathcal{R}}_a\rho(r,t)\partial_r\rho(r,t)=0,
\end{align}
where the one-dimensional transform $\widetilde{\mathcal{R}}_a$ is given by
\begin{align*}
\widetilde{\mathcal{R}}_a\rho(r,t)
&=\frac{\Gamma(\frac{n-1}{2})}{2\pi^{\frac{n+1}{2}}}\int\limits_0^\infty\partial_\zeta\rho(\zeta,t)\zeta^{n-1}
\int\limits_{\mathbb{S}^{n-1}}\Big(\frac{y_1}{(|re_1-\zeta y|^2+a^2)^{\frac{n-1}{2}}}-\frac{y_1}{|re_1-\zeta y|^{n-1}}\Big)d\sigma(y)d\zeta.
\end{align*}
Note that \eqref{radial-equation} is a one-dimensional transport equation. Following the flow map arguments in \cite{[Kiselev-Sarsam]}, we can derive that
\begin{align}\label{derivative-formula}
\partial_r\rho(r,t)=\rho'_0(\Phi^{-1}_t(r))e^{-g\int_0^t\widetilde{\mathcal{R}}_a(\partial_r\rho)(\Phi_s\circ\Phi^{-1}_t(r),s)ds},
\end{align}
where the flow map $\Phi_t(r)$ is defined by
\begin{equation*}
\left\{
\begin{array}{ll}
\frac{d}{dt}\Phi_t(r)=g\widetilde{\mathcal{R}}_a\rho(\Phi_t(r),r),\\
\Phi_0(r)=r,
\end{array}
\right.
\end{equation*}
for each fixed $r\in[0,\infty)$. By the assumption $\rho'_0\geq0$, \eqref{derivative-formula} shows that $\rho$ is radially non-decreasing.
The proof of Lemma \ref{radial-property-preserved} is then finished. \hfill\hfill$\square$\vskip12pt

\subsection{A positive lower bound for the nonlinear term}
Next we derive a positive lower bound for the nonlinear term of \eqref{M-KS}, which is vital for the proof of finite time blow-up. For this purpose, we need a pointwise inequality for the transform $\mathcal{R}_a$. The similar inequality for the one-dimensional transform $\mathcal{H}_a$ was established in \cite{[Liu-Zhang]}.
\begin{lemma}\label{lower-bound-Ra}
 Fix $n\geq2$ and $a>0$. Let $f:\mathbb{R}^n\rightarrow\mathbb{R}$ be a radial, non-decreasing and continuously differentiable function with $\nabla f\in L^1(\mathbb{R}^n)\cap L^\infty(\mathbb{R}^n)$. Then,  for any $x\neq0$, we have
\begin{eqnarray*}
-\mathcal{R}_a f(x)\cdot\frac{x}{|x|}
\geq\frac{nB(\frac12,\frac{n+1}{2})}{2^{n+1}\pi}\frac{1}{|x|^n}\Big(1-\frac{2^{n+1}|x|^{n+1}}{(4|x|^2+a^2)^{\frac{n+1}{2}}}\Big)
\int\limits_0^{|x|} (f(|x|)-f(\varrho))\varrho^{n-1}
d\varrho.
\end{eqnarray*}
\end{lemma}
\begin{remark}
For a radial and non-increasing Schwartz function $f:\mathbb{R}^n\rightarrow\mathbb{R}$, we have the inequality
\begin{eqnarray*}
\mathcal{R}_a f(x)\cdot\frac{x}{|x|}
\geq\frac{nB(\frac12,\frac{n+1}{2})}{2^{n+1}\pi}\frac{1}{|x|^n}\Big(1-\frac{2^{n+1}|x|^{n+1}}{(4|x|^2+a^2)^{\frac{n+1}{2}}}\Big)
\int\limits_0^{|x|} (f(\varrho)-f(|x|))\varrho^{n-1}
d\varrho.
\end{eqnarray*}
\end{remark}
\textbf{Proof of Lemma \ref{lower-bound-Ra}.} By \eqref{definition-Ra} and \eqref{Ra-kernel}, we note that
\begin{align*}
\mathcal{R}_af(x)
&=\mathcal{R}f(x)-\frac{\Gamma(\frac{n+1}{2})}{\pi^{\frac{n+1}{2}}}\int\limits_{\mathbb{R}^n}\frac{x-y}{(|x-y|^2+a^2)^{\frac{n+1}{2}}}f(y)dy\\
&=-\Lambda^{-1}\nabla f(x)-\frac{\Gamma(\frac{n+1}{2})}{\pi^{\frac{n+1}{2}}}\int\limits_{\mathbb{R}^n}\frac{x-y}{(|x-y|^2+a^2)^{\frac{n+1}{2}}}f(y)dy.
\end{align*}
Then, by the integral representation of the Riesz potential, integration by parts and the radial assumption on $f$, we have
\begin{align*}
\mathcal{R}_af(x)
&=-\frac{\Gamma(\frac{n-1}{2})}{2\pi^{\frac{n+1}{2}}}\int\limits_{\mathbb{R}^{n}}\frac{\nabla f(y)}{|x-y|^{n-1}}dy+\frac{\Gamma(\frac{n+1}{2})}{(n-1)\pi^{\frac{n+1}{2}}}\int\limits_{\mathbb{R}^{n}}\frac{\nabla f(y)}{(|x-y|^2+a^2)^{\frac{n-1}{2}}}dy\\
&=-\frac{\Gamma(\frac{n-1}{2})}{2\pi^{\frac{n+1}{2}}}\int\limits_{\mathbb{R}^{n}}\Big(\frac{1}{|x-y|^{n-1}}-\frac{1}{(|x-y|^2+a^2)^{\frac{n-1}{2}}}\Big)\nabla f(y)dy\\
&=-\frac{\Gamma(\frac{n-1}{2})}{2\pi^{\frac{n+1}{2}}}\int\limits_{\mathbb{R}^{n}}\Big(\frac{1}{|x-y|^{n-1}}-\frac{1}{(|x-y|^2+a^2)^{\frac{n-1}{2}}}\Big)f'(|y|)\frac{y}{|y|}dy\\
&=-\frac{\Gamma(\frac{n-1}{2})}{2\pi^{\frac{n+1}{2}}}\int\limits_0^\infty f'(\varrho)\varrho^{n-1}\int\limits_{\mathbb{S}^{n-1}}\Big(\frac{z}{|x-\varrho z|^{n-1}}-\frac{z}{(|x-\varrho z|^2+a^2)^{\frac{n-1}{2}}}\Big)d\sigma(z)d\varrho,
\end{align*}
which follows from the rotational transform and a change of variables  that, for any $x\neq0$,
\begin{align*}
-\mathcal{R}_af(x)\cdot\frac{x}{|x|}
=\frac{\Gamma(\frac{n-1}{2})}{2\pi^{\frac{n+1}{2}}}\int\limits_0^\infty f'(\varrho)\varrho^{n-1}
\int\limits_{\mathbb{S}^{n-1}}\Big(\frac{z_1}{||x|e_1-\varrho z|^{n-1}}-\frac{z_1}{(||x|e_1-\varrho z|^2+a^2)^{\frac{n-1}{2}}}\Big)d\sigma(z)d\varrho.
\end{align*}
By a change of variables formula (see e.g., pp. 592 of \cite{[Grafakos]}), integration by parts and the mean value theorem, we observe that, for $\varrho\neq|x|$,
\begin{align*}
&\int\limits_{\mathbb{S}^{n-1}}\Big(\frac{z_1}{||x|e_1-\varrho z|^{n-1}}-\frac{z_1}{(||x|e_1-\varrho z|^2+a^2)^{\frac{n-1}{2}}}\Big)d\sigma(z)\\
&=\int\limits_{-1}^1\int\limits_{\sqrt{1-s^2}\mathbb{S}^{n-2}}\Big(\frac{s}{((|x|-\varrho s)^2+\varrho^2|z|^2)^{\frac{n-1}{2}}}-\frac{s}{((|x|-\varrho s)^2+\varrho^2|z|^2+a^2)^{\frac{n-1}{2}}}\Big)d\sigma(z)
\frac{ds}{\sqrt{1-s^2}}\\
&=\omega_{n-2}\int\limits_{-1}^1\Big(\frac{s(1-s^2)^{\frac{n-3}{2}}}{(|x|^2-2|x|\varrho s+\varrho^2)^{\frac{n-1}{2}}}-\frac{s(1-s^2)^{\frac{n-3}{2}}}{(|x|^2-2|x|\varrho s+\varrho^2+a^2)^{\frac{n-1}{2}}}\Big)ds\\
&=\omega_{n-2}\int\limits_0^\pi\Big(\frac{\cos\mu\sin^{n-2}\mu }{(|x|^2-2|x|\varrho\cos\mu+\varrho^2)^{\frac{n-1}{2}}}-\frac{\cos\mu\sin^{n-2}\mu }{(|x|^2-2|x|\varrho\cos\mu+\varrho^2+a^2)^{\frac{n-1}{2}}}\Big)d\mu\\
&=\omega_{n-2}\varrho|x|\int\limits_0^\pi\Big(\frac{\sin^n\mu }{(|x|^2-2|x|\varrho\cos\mu+\varrho^2)^{\frac{n+1}{2}}}-\frac{\sin^n\mu }{(|x|^2-2|x|\varrho\cos\mu+\varrho^2+a^2)^{\frac{n+1}{2}}}\Big)d\mu\\
&=\omega_{n-2}\frac{\varrho}{|x|^{n}}\int\limits_0^\pi\Big(\frac{\sin^n\mu }{(1-2\frac{\varrho}{|x|}\cos\mu+\frac{\varrho^2}{|x|^2})^{\frac{n+1}{2}}}-\frac{\sin^n\mu }{(1-2\frac{\varrho}{|x|}\cos\mu+\frac{\varrho^2}{|x|^2}+\frac{a^2}{|x|^2})^{\frac{n+1}{2}}}\Big)d\mu\\
&=\omega_{n-2}\frac{n+1}{2}\frac{\varrho}{|x|^{n}}\frac{a^2}{|x|^2}\int\limits_0^\pi \int\limits_0^1\frac{\sin^n\mu}{(1-2\frac{\varrho}{|x|}\cos\mu+\frac{\varrho^2}{|x|^2}+\tau\frac{a^2}{|x|^2})^{\frac{n+3}{2}}} d\tau d\mu,
\end{align*}
where $\omega_{n-2}=\frac{2\pi^{\frac{n-1}{2}}}{\Gamma(\frac{n-1}{2})}$ is the surface area of $\mathbb{S}^{n-2}$.
Thus, we obtain that
\begin{align*}
-\mathcal{R}_af(x)\cdot\frac{x}{|x|}
&=\frac{n+1}{2\pi}\frac{a^2}{|x|^{n+2}}\int\limits_0^\infty f'(\varrho)\varrho^{n}
\int\limits_0^\pi \int\limits_0^1\frac{\sin^n\mu}{(1-2\frac{\varrho}{|x|}\cos\mu+\frac{\varrho^2}{|x|^2}+\tau\frac{a^2}{|x|^2})^{\frac{n+3}{2}}} d\tau d\mu d\varrho\\
&\geq\frac{n+1}{2\pi}\frac{a^2}{|x|^{n+2}}\int\limits_0^{|x|} f'(\varrho)\varrho^{n}
\int\limits_0^\pi \int\limits_0^1\frac{\sin^n\mu}{(1-2\frac{\varrho}{|x|}\cos\mu+\frac{\varrho^2}{|x|^2}+\tau\frac{a^2}{|x|^2})^{\frac{n+3}{2}}} d\tau d\mu d\varrho\\
&\geq\frac{n+1}{2\pi}\frac{a^2}{|x|^{n+2}}\int\limits_0^{|x|} f'(\varrho)\varrho^{n}
\int\limits_0^\pi \int\limits_0^1\frac{\sin^n\mu}{(4+\tau\frac{a^2}{|x|^2})^{\frac{n+3}{2}}} d\tau d\mu d\varrho\\
&=\frac{B(\frac12,\frac{n+1}{2})}{2^{n+1}\pi}\frac{1}{|x|^n}\Big(1-\frac{2^{n+1}|x|^{n+1}}{(4|x|^2+a^2)^{\frac{n+1}{2}}}\Big)\int\limits_0^{|x|} f'(\varrho)\varrho^{n}d\varrho\\
&=\frac{nB(\frac12,\frac{n+1}{2})}{2^{n+1}\pi}\frac{1}{|x|^n}\Big(1-\frac{2^{n+1}|x|^{n+1}}{(4|x|^2+a^2)^{\frac{n+1}{2}}}\Big)
\int\limits_0^{|x|} (f(|x|)-f(\varrho))\varrho^{n-1}
d\varrho,
\end{align*}
which is the desired lower bound. \hfill\hfill$\square$\vskip12pt
\begin{proposition}\label{n-nonlinear-inequality-monotone}
Fix $n\geq2$ and $a>0$. Let $-1<\delta<1$ and $f:\mathbb{R}^n\rightarrow\mathbb{R}$ be a radial, non-decreasing and continuously differentiable function with $\nabla f\in L^1(\mathbb{R}^n)\cap L^\infty(\mathbb{R}^n)$.
Then
\begin{eqnarray}\label{positive-lower-bound}
-\int\limits_{\mathbb{R}^n}\frac{\mathcal{R}_a f(x)\cdot\nabla f(x)}{|x|^{n+\delta}}dx
\geq C_{n,\delta}
\int\limits_{\mathbb{R}^n}\frac{(f(x)-f(0))^2}{|x|^{n+1+\delta}}\Big(1-\frac{2^{n+1}|x|^{n+1}}{(4|x|^2+a^2)^{\frac{n+1}{2}}}\Big)dx,
\end{eqnarray}
where
\begin{align*}
C_{n,\delta}=\frac{(\sqrt{n+1+\delta}-\sqrt{n})^2}{2^{n+2}\pi}B\Big(\frac12,\frac{n+1}{2}\Big).
\end{align*}
\end{proposition}
\textbf{Proof.} By Lemma \ref{lower-bound-Ra} and the monotonicity of $f$, we have that
\begin{eqnarray*}
\begin{split}
-\mathcal{R}_af(x)\cdot\nabla f(x)
&=-\mathcal{R}_af(x)\cdot\frac{x}{|x|} f'(|x|)\\
&\geq\frac{nB(\frac12,\frac{n+1}{2})}{2^{n+1}\pi}\frac{f'(|x|)}{|x|^n}\Big(1-\frac{2^{n+1}|x|^{n+1}}{(4|x|^2+a^2)^{\frac{n+1}{2}}}\Big)\int\limits_0^{|x|} (f(|x|)-f(\varrho))\varrho^{n-1}
d\varrho,
\end{split}
\end{eqnarray*}
which implies that
\begin{align*}
-\int\limits_{\mathbb{R}^n}\frac{\mathcal{R}_a f(x)\cdot\nabla f(x)}{|x|^{n+\delta}}dx
&\geq\frac{nB(\frac12,\frac{n+1}{2})}{2^{n+1}\pi}\int\limits_{\mathbb{R}^n}\frac{f'(|x|)}{|x|^{2n+\delta}}\Big(1-\frac{2^{n+1}|x|^{n+1}}{(4|x|^2+a^2)^{\frac{n+1}{2}}}\Big)\int\limits_0^{|x|} (f(|x|)-f(\varrho))\varrho^{n-1}
d\varrho dx\\
&=\frac{nB(\frac12,\frac{n+1}{2})}{2^{n+1}\pi}\omega_{n-1}\int\limits_0^\infty\frac{f'(r)}{r^{n+1+\delta}}\Big(1-\frac{2^{n+1}r^{n+1}}{(4r^2+a^2)^{\frac{n+1}{2}}}\Big)\int\limits_0^{r} (f(r)-f(\varrho))\varrho^{n-1}
d\varrho dr.
\end{align*}
Furthermore, by Fubini's theorem and integration by parts, we derive that
\begin{align*}
&\int\limits_0^\infty\frac{f'(r)}{r^{n+1+\delta}}\Big(1-\frac{2^{n+1}r^{n+1}}{(4r^2+a^2)^{\frac{n+1}{2}}}\Big)\int\limits_0^{r} (f(r)-f(\varrho))\varrho^{n-1}
d\varrho dr\\
&=\frac12\int\limits_0^\infty\varrho^{n-1}\int\limits_\varrho^\infty\frac{1}{r^{n+1+\delta}}\Big(1-\frac{2^{n+1}r^{n+1}}{(4r^2+a^2)^{\frac{n+1}{2}}}\Big)\frac{\partial}{\partial r}\Big((f(r)-f(\varrho))^2\Big) drd\varrho\\
&=\frac{n+1+\delta}{2}\iint\limits_{r\geq\varrho>0}\frac{(f(r)-f(\varrho))^2}{r^{n+2+\delta}}\Big(1-\frac{2^{n+1}r^{n+1}}{(4r^2+a^2)^{\frac{n+1}{2}}}\Big)\varrho^{n-1}d\varrho dr\\
&\ \ \ \
+2^n(n+1)\iint\limits_{r\geq\varrho>0}\frac{a^2(f(r)-f(\varrho))^2}{r^{1+\delta}(4r^2+a^2)^{\frac{n+3}{2}}}\varrho^{n-1}d\varrho dr\\
&\geq\frac{n+1+\delta}{2}\iint\limits_{r\geq\varrho>0}\frac{(f(r)-f(\varrho))^2}{r^{n+2+\delta}}\Big(1-\frac{2^{n+1}r^{n+1}}{(4r^2+a^2)^{\frac{n+1}{2}}}\Big)\varrho^{n-1}d\varrho dr.
\end{align*}
Now using the elementary inequality
\begin{eqnarray}\label{young}
(b_1-b_2)^2\geq(1-\alpha)b^2_1+\Big(1-\frac{1}{\alpha}\Big)b^2_2
\end{eqnarray}
for any $b_1,b_2\in\mathbb{R}$ and any $0<\alpha<1$, we obtain
\begin{align*}
&\iint\limits_{r\geq\varrho>0}\frac{(f(r)-f(\varrho))^2}{r^{n+2+\delta}}\Big(1-\frac{2^{n+1}r^{n+1}}{(4r^2+a^2)^{\frac{n+1}{2}}}\Big)\varrho^{n-1}d\varrho dr\\
&\geq(1-\alpha)\iint\limits_{r\geq\varrho>0}\frac{(f(r)-f(0))^2}{r^{n+2+\delta}}\Big(1-\frac{2^{n+1}r^{n+1}}{(4r^2+a^2)^{\frac{n+1}{2}}}\Big)\varrho^{n-1}d\varrho dr\\
&\ \ \ \ \
+\Big(1-\frac{1}{\alpha}\Big)\iint\limits_{r\geq\varrho>0}\frac{(f(\varrho)-f(0))^2}{r^{n+2+\delta}}\Big(1-\frac{2^{n+1}r^{n+1}}{(4r^2+a^2)^{\frac{n+1}{2}}}\Big)\varrho^{n-1}d\varrho dr\\
&\geq(1-\alpha)\iint\limits_{r\geq\varrho>0}\frac{(f(r)-f(0))^2}{r^{n+2+\delta}}\Big(1-\frac{2^{n+1}r^{n+1}}{(4r^2+a^2)^{\frac{n+1}{2}}}\Big)\varrho^{n-1}d\varrho dr\\
&\ \ \ \ \
+\Big(1-\frac{1}{\alpha}\Big)\iint\limits_{r\geq\varrho>0}\frac{(f(\varrho)-f(0))^2}{r^{n+2+\delta}}\Big(1-\frac{2^{n+1}\varrho^{n+1}}{(4\varrho^2+a^2)^{\frac{n+1}{2}}}\Big)\varrho^{n-1}d\varrho dr\\
&=\frac{1-\alpha}{n}\int\limits^\infty_0\frac{(f(r)-f(0))^2}{r^{2+\delta}}\Big(1-\frac{2^{n+1}r^{n+1}}{(4r^2+a^2)^{\frac{n+1}{2}}}\Big)dr\\
&\ \ \ \ \
+\frac{1-\frac{1}{\alpha}}{n+1+\delta}\int\limits^\infty_0\frac{(f(\varrho)-f(0))^2}{\varrho^{2+\delta}}\Big(1-\frac{2^{n+1}\varrho^{n+1}}{(4\varrho^2+a^2)^{\frac{n+1}{2}}}\Big)d\varrho\\
&=\frac{2n+1+\delta}{n(n+1+\delta)}\int\limits^\infty_0\frac{(f(r)-f(0))^2}{r^{2+\delta}}\Big(1-\frac{2^{n+1}r^{n+1}}{(4r^2+a^2)^{\frac{n+1}{2}}}\Big)dr\\
&\ \ \
-\Big(\frac{\alpha}{n}+\frac{1}{(n+1+\delta)\alpha}\Big)\int\limits^\infty_0\frac{(f(r)-f(0))^2}{r^{2+\delta}}\Big(1-\frac{2^{n+1}r^{n+1}}{(4r^2+a^2)^{\frac{n+1}{2}}}\Big)dr,
\end{align*}
where in the second inequality we have used the monotonicity of the weight $1-\frac{2^{n+1}r^{n+1}}{(4r^2+a^2)^{\frac{n+1}{2}}}$ in $r$ and the fact that $1-\frac{1}{\alpha}<0$.
It follows from the choice of $\alpha=\sqrt{\frac{n}{n+1+\delta}}\in(0,1)$ that
\begin{align*}
&\iint\limits_{r\geq\varrho>0}\frac{(f(r)-f(\varrho))^2}{r^{n+2+\delta}}\Big(1-\frac{2^{n+1}r^{n+1}}{(4r^2+a^2)^{\frac{n+1}{2}}}\Big)\Big)\varrho^{n-1}d\varrho dr\\
&\geq\Big(\frac{2n+1+\delta}{n(n+1+\delta)}-\frac{2}{\sqrt{n(n+1+\delta)}}\Big)\int\limits^\infty_0\frac{(f(r)-f(0))^2}{r^{2+\delta}}\Big(1-\frac{2^{n+1}r^{n+1}}{(4r^2+a^2)^{\frac{n+1}{2}}}\Big)dr\\
&=\frac{(\sqrt{n+1+\delta}-\sqrt{n})^2}{n(n+1+\delta)}\int\limits^\infty_0\frac{(f(r)-f(0))^2}{r^{2+\delta}}\Big(1-\frac{2^{n+1}r^{n+1}}{(4r^2+a^2)^{\frac{n+1}{2}}}\Big)dr.
\end{align*}
Altogether, we can obtain that
\begin{align*}
-\int\limits_{\mathbb{R}^n}\frac{\mathcal{R}_a f(x)\cdot\nabla f(x)}{|x|^{n+\delta}}dx
\geq\frac{(\sqrt{n+1+\delta}-\sqrt{n})^2}{2^{n+2}\pi}B(\frac12,\frac{n+1}{2})\int\limits_{\mathbb{R}^n}\frac{(f(x)-f(0))^2}{|x|^{n+1+\delta}}\Big(1-\frac{2^{n+1}|x|^{n+1}}{(4|x|^2+a^2)^{\frac{n+1}{2}}}\Big)dx,
\end{align*}
which concludes the proof of Proposition \ref{n-nonlinear-inequality-monotone}. \hfill\hfill$\square$\vskip12pt
\begin{remark}
For a radial and non-increasing Schwartz function $f:\mathbb{R}^n\rightarrow\mathbb{R}$, we also have the bilinear inequality \eqref{positive-lower-bound}.
\end{remark}
\subsection{Proof of Theorem \ref{singularity-formation-monotone} (blow-up of solutions)}
With the help of Proposition \ref{n-nonlinear-inequality-monotone}, we are ready to prove Theorem \ref{singularity-formation-monotone}.

\textbf{Proof of Theorem \ref{singularity-formation-monotone}.}
Suppose that the initial data $\rho_0\in C^\infty_c(\mathbb{R}^n)$ is radial and non-decreasing with $\rho_0(0)<0$. Moreover, suppose that $\rho_0$ is supported in a ball $B_L(0)$ for some $L>0$.
Let $\rho(x,t)$ denote the corresponding unique local-in-time solution to \eqref{M-KS}.
The Beale-Kato-Majda type criterion in Proposition \ref{Beale-Kato-Majda-type} reduces to the proof to show that the class of initial data always leads to finite time blow-up in some way.

For this purpose, assume for contradiction that $\rho$ exists for all time. By Lemma \ref{radial-property-preserved}, $\rho$ is radial non-decreasing on $\mathbb{R}^n$ for all time $t>0$.
Then,  it is not difficult to see from the representation formula of $-\mathcal{R}_a\rho(x,t)\cdot\frac{x}{|x|}$ in the proof of Lemma \ref{lower-bound-Ra} that $\mathcal{R}_a\rho(x,t)\cdot(-\frac{x}{|x|})\geq0$ for any $x\neq0$ and $t>0$. This shows that the velocity field $g\mathcal{R}_a\rho(x,t)$ points inwards along the boundary $\partial B_L(0)$. Therefore, for any $t>0$ the support of $\rho(\cdot,t)$ lies inside $B_L(0)$.

We proceed to denote
\begin{eqnarray*}
I(t):=\int\limits_{B_L(0)}\frac{\rho(x,t)-\rho(0,t)}{|x|^{n+\delta}}dx,
\end{eqnarray*}
where $\delta\in(0,1)$ is arbitrarily fixed.
By H\"{o}lder's inequality, we obtain that
\begin{align*}
|I(t)|
\leq \|\nabla\rho\|_{L^\infty}\int_{|x|\leq L}\frac{dx}{|x|^{n-1+\delta}}
=\frac{\omega_{n-1}L^{1-\delta}}{1-\delta}\|\nabla\rho\|_{L^\infty}<+\infty,
\end{align*}
which shows that $I(t)$ is finite for all time $t>0$.

Next we prove that $I(t)$ will blow up in finite time, leading a contradiction.
To this end, by Lemma \ref{radial-property-preserved}, we know that the velocity at the origin is $0$, that is,
\begin{eqnarray*}
\mathcal{R}_a\rho(0,t)=\frac{\Gamma(\frac{n+1}{2})}{\pi^{\frac{n+1}{2}}}P.V.\int\limits_{\mathbb{R}^n}\Big(\frac{y}{(|y|^2+a^2)^{\frac{n+1}{2}}}-\frac{y}{|y|^{n+1}}\Big)\rho(|y|,t)dy=0,
\end{eqnarray*}
which yields that $\rho(0,t)=\rho_0(0)$ for any $t>0$. Then, we multiply the equation \eqref{M-KS} by the weight $|x|^{-n-\delta}$ and integrate in $x\in B_L(0)$. By Proposition \ref{n-nonlinear-inequality-monotone}, we compute
\begin{align*}
\frac{d}{dt}I(t)
&=\int\limits_{B_L(0)}\frac{\partial_t\rho(x,t)}{|x|^{n+\delta}}dx=-g\int\limits_{B_L(0)}\frac{\mathcal{R}_a \rho(x,t)\cdot\nabla \rho(x,t)}{|x|^{n+\delta}}dx
=-g\int\limits_{\mathbb{R}^n}\frac{\mathcal{R}_a \rho(x,t)\cdot\nabla \rho(x,t)}{|x|^{n+\delta}}dx\\
&\geq gC_{n,\delta}
\int\limits_{\mathbb{R}^n}\frac{(\rho(x,t)-\rho(0,t))^2}{|x|^{n+1+\delta}}\Big(1-\frac{2^{n+1}|x|^{n+1}}{(4|x|^2+a^2)^{\frac{n+1}{2}}}\Big)dx\\
&\geq gC_{n,\delta}
\int\limits_{B_L(0)}\frac{(\rho(x,t)-\rho(0,t))^2}{|x|^{n+1+\delta}}\Big(1-\frac{2^{n+1}|x|^{n+1}}{(4|x|^2+a^2)^{\frac{n+1}{2}}}\Big)dx\\
&\geq gC_{n,\delta}
\Big(1-\frac{2^{n+1}L^{n+1}}{(4L^2+a^2)^{\frac{n+1}{2}}}\Big)\int\limits_{B_L(0)}\frac{(\rho(x,t)-\rho(0,t))^2}{|x|^{n+1+\delta}}dx\\
&\geq gC_{n,\delta}
\Big(1-\frac{2^{n+1}L^{n+1}}{(4L^2+a^2)^{\frac{n+1}{2}}}\Big)\Big(\int\limits_{B_L(0)}\frac{dx}{|x|^{n-1+\delta}}\Big)^{-1}(I(t))^2\\
&=\frac{g(1-\delta)C_{n,\delta}}{\omega_{n-1}L^{1-\delta}}
\Big(1-\frac{2^{n+1}L^{n+1}}{(4L^2+a^2)^{\frac{n+1}{2}}}\Big)(I(t))^2,
\end{align*}
where in the last inequality we have used the fact that $0<\delta<1$ and the Cauchy-Schwarz inequality.
Since
\begin{eqnarray*}
I(0)=\int\limits_{B_L(0)}\frac{\rho_0(x)-\rho_0(0)}{|x|^{n+\delta}}dx>0,
\end{eqnarray*}
it follows that $I(t)$ must blow up in finite time.
We then have proved the theorem.   \hfill\hfill$\square$\vskip12pt
\subsection*{Acknowledgment}
\addcontentsline{toc}{section}{Acknowledgments} \hskip\parindent
W. Zhang was supported by the National Natural Science Foundation of China (NNSFC) (No. 12501301).





\end{document}